\documentclass{article}

\newtheorem{theorem}{Theorem}

\newtheorem{conjecture}[theorem]{Conjecture}

\input{tcilatex}

\begin{document}

\begin{center}
\bigskip

{\Huge A generalization of the Collatz problem and conjecture}

\bigskip

{\Large M. Bruschi}

Dipartimento di Fisica, Universit\`{a} di Roma "La Sapienza", Rome, Italy

Istituto Nazionale di Fisica Nucleare, Sezione di Roma

$\emph{mario.bruschi@roma1.infn.it}$

\bigskip

\textit{Summary}
\end{center}

We introduce an infinite set of integer mappings that generalize the
well-known Collatz-Ulam mapping and we conjecture that an infinite subset of
these mappings feature the remarkable property of the Collatz conjecture,
namely that they converge to unity irrespective of which positive integer is
chosen initially.

\newpage

\section{Introduction}

The Collatz problem is so well known that we refer for formulation,
references and bibliography to the web \cite{wiki},\cite{wolf},\cite{lag}.

Let us quote from ref. \cite{wiki}:

"The Collatz conjecture is an unsolved conjecture in mathematics. It is
named after Lothar Collatz, who first proposed it in 1937. The conjecture is
also known as the $3n+1$ conjecture, as the Ulam conjecture (after Stanislaw
Ulam), or as the Syracuse problem; the sequence of numbers involved is
referred to as the hailstone sequence or hailstone numbers, or as wondrous
numbers per \textit{G\"{o}del, Escher, Bach \cite{bge}.}

We take any number $n$. If $n$ is even, we halve it ($n/2$), else we do
"triple plus one" and get $3n+1$. The conjecture is that for all numbers
this process converges to $1$. Hence it has been called 'Half Or Triple Plus
One', sometimes called HOTPO.

Paul Erd\H{o}s said about the Collatz conjecture: 'Mathematics is not yet
ready for such problems.' He offered \$500 for its solution. "\medskip

In this paper we present a generalization of the Collatz map, and prove the
corresponding conjecture for a  set of initial values; this proof has no
relevance for the original Collatz conjecture, although one might hope that
it provide some hint for solving that problem.

\subsection{The Collatz problem}

Consider the following operation on an arbitrary positive integer:

* If the number is even, divide it by two.

* If the number is odd, triple it and add one.

For example, if this operation is performed on 3, the result is 10; if it is
performed on 28, the result is 14.

In modular arithmetic notation, define the function $f(n)$ as follows ($n$
is a positive integer):

\begin{eqnarray}
f(n) &=&n/2~\ \ \ \ \ \ \ \ \ if\ \ n\equiv 0\ \ (\func{mod}\ 2)~,  \nonumber
\\
f(n) &=&3n+1~\ \ \ \ ~if\ \ n\equiv 1\ \ (\func{mod}\ 2)~.  \label{cmap}
\end{eqnarray}

Now, form a sequence $S_{k}$ by performing this operation repeatedly,
beginning with any positive integer $n$, and taking the result at each step
as the input for the next.

In notation ($n$ is a positive integer):%
\begin{equation}
S_{0}=n~;\ ~~~~~S_{k}=f(S_{k-1})~,~k=1,2,...~.  \label{SeqC}
\end{equation}

The Collatz conjecture can then be formulated as follows:

\begin{conjecture}
This process will eventually reach the number $1$, irrespective of which
positive integer $n$ is chosen initially.
\end{conjecture}

Note that, once the value $n=1$ is reached, the sequence (\ref{SeqC})
reduces to the limit cycle$~$yielding sequentially the values $1,~4,~2,~1$.

This Conjecture is simple to state but so difficult to prove that the issue
of its validity is still open.

Many generalizations of the Collatz function (\ref{cmap}) have been
introduced: but, to the best of our knowledge, the generalizations that we
propose in this paper are new.

In the following Section we introduce an infinite set of generalized Collatz
maps (involving two additional integer parameters $b,m)$; in our opinion
they qualify as "correct" generalizations of the original Collatz map,
inasmuch as, for most of these maps, the characteristic property of the
original Collatz Conjecture seems to hold. We can also prove, albeit only
for a limited set of values of the starting integer $n$, that these
generalized maps do indeed possess the Collatz property ("convergence to
unity"),

\subsection{Generalizations}

Define the function $f(n,b,m)$ as follows ($n,m=1,2,..;b=2,3,...$):

\begin{eqnarray}
f(n,b,m) &=&n/b~\ \ \ \ \ \ \ \ \ if\ \ n\equiv 0\ \ (\func{mod}\ b)~,
\label{rule1} \\
f(n,b,m) &=&(b^{m}~+1)n+b^{m}~-(n~\func{mod}\ b^{m}~)\ \ \ \ ~if\ \ n\neq 0\
\ (\func{mod}\ b)~.  \label{rule2}
\end{eqnarray}

Now, form a sequence $S_{k}$ by performing this operation repeatedly,
beginning with any positive integer $n$, and taking the result at each step
as the input for the next step .

In notation ($n$ is a positive integer):

\begin{equation}
S_{0}=n;~~~~~S_{k}=f(S_{k-1},b,m)\ ,~k=1,2,...  \label{SeqB}
\end{equation}%
\medskip

\textbf{Remark} \ \textit{Clearly for }$b=2,m=1$\textit{\ one gets the
original Collatz map.\smallskip }

We can easily prove the limited result given by the following proposition:

\textbf{Proposition} \textit{\ If }$S_{0}<b^{m}$\textit{\ then the sequence (%
\ref{SeqB}) shall eventually reach the number }$1$\textit{\ (for any }$%
b>1,m>0$\textit{)}

\bigskip \textbf{Proof}

Let us write the integer $S_{0}<b^{m}$ in base $b:$%
\begin{equation}
S_{0}=a_{0}+a_{1}b+a_{2}b^{2}+..+a_{m-1}b^{m-1}  \label{s0}
\end{equation}%
where of course $a_{i}<b,i=0,1,..,m-1$ and let us assume without loss of
generality that $a_{0}>0$ (otherwise we apply the rule (\ref{rule1})) and
for convenience let us put 
\begin{equation}
a_{0}=b-j,~0<j<b.  \label{j}
\end{equation}

. Then the rule (\ref{rule2}) yields%
\begin{eqnarray}
S_{1} &=&(b^{m}~+1)S_{0}+b^{m}~-(S_{0}~\func{mod}\ b^{m}~)\  \\
&=&b^{m}~\left( S_{0}+1\right) +S_{0}-(S_{0}~\func{mod}\ b^{m}~).
\end{eqnarray}%
Now, noting that (see (\ref{s0})) 
\begin{equation}
S_{0}~\func{mod}\ b^{m}=S_{\text{ }0}
\end{equation}%
we have%
\begin{equation}
S_{1}=b^{m}\left[ \left( a_{0}+1\right) +a_{1}b+a_{2}b^{2}+..+a_{m-1}b^{m-1}%
\right] .
\end{equation}%
Let us assume, just for the sake of clearness in the exposition, that \
\bigskip $j\neq 1$ so that $\left( a_{0}+1\right) <b$ (see (\ref{j})). Then
we have to apply $m$ times the rule (\ref{rule1}), reaching%
\begin{equation}
S_{m+1}=\left( a_{0}+1\right) +a_{1}b+a_{2}b^{2}+..+a_{m-1}b^{m-1}.
\end{equation}%
Repeating the above procedure, one gets (see (\ref{j})) 
\begin{eqnarray}
S_{j\cdot m+j} &=&\left( a_{0}+j\right) +a_{1}b+a_{2}b^{2}+..+a_{m-1}b^{m-1}
\\
&=&\left( b-j+j\right) +a_{1}b+a_{2}b^{2}+..+a_{m-1}b^{m-1} \\
&=&\left( a_{1}+1\right) b+a_{2}b^{2}+..+a_{m-1}b^{m-1}.
\end{eqnarray}%
Then the rule (\ref{rule1}) yields%
\begin{equation}
S_{j\cdot m+j+1}=\left( a_{1}+1\right) +a_{2}b^{2-1}+..+a_{m-1}b^{m-2}.
\end{equation}

It is now clear that eventually the sequence reach the number $%
1.(QED)\medskip $

\textbf{Remark} \ \ \textit{The above proof makes also clear how \ to
compute \ the\ stopping time (namely the number of steps necessary to reach\
the number }$1$\textit{)\ for the sequences starting from }$%
S_{0}<b^{m}.\smallskip $

\textbf{Remark \ \ }\textit{It is also interesting to see the limit cycles,
namely the cyclic sequences obtained starting from }$S_{0}=1$\textit{\ for
different values of \ }$b,m$\textit{. They are:}%
\begin{equation}
1,(2b^{m},2b^{m-1},..,2b),2,(3b^{m},3b^{m-1},..,3b),3,......,b-1,..,b,1.
\label{LC}
\end{equation}

\textit{\ Just as an example we report the limit cycle for }$b=5,m=3:$ 
\begin{equation}
1,250,50,10,2,375,75,15,3,500,100,20,4,625,125,25,5,1.
\end{equation}

\textbf{\ \ }

The above proposition covers a limited set of starting points for the
sequence (a very small set indeed for small values of \ $b,m$\ but a rapidly
increasing one for increasing values of $b,m$). Of course this set can be
trivially extended to an infinite set considering all the values $%
S_{0}=s~b^{N},~s<b^{m},~N=0,1,2,..$ For all these initial values the
sequence (\ref{SeqB}) obviously converges to the number $1$ and thus the
limit cycle (\ref{LC}) due to the above Proposition and to the rule (\ref%
{rule1}).\smallskip\ Let us call '\textit{trivial}' all these values of $%
S_{0}$ (and "\textit{non-trivial}\emph{" }the values outside this set). 

An obvious question arises: since our mapping is a natural generalization of
the Collatz one and considering the result of the above Proposition, may we
hope that a generalized Collatz conjecture hold true, namely that for any $%
b,m,S_{0}$ the limit cycle (\ref{LC}) is eventually reached?\smallskip
\medskip 

Having failed (up to now) to extend the result of the above Proposition, we
had to resort to a massive computer investigation which of course proves
nothing... moreover it is very time-expensive even for relatively small
values of $b,m$ if one tries to test a substantial number (say $10^{6}$-$%
10^{9}$) of \textit{non-trivial }values of $S_{0}$ (consider for instance
that with $b=10,m=9$ for the first \textit{non-trivial }value of $S_{0}$,
namely $S_{0}=10^{9}+1,$ the \emph{stopping time} is $\ 5000000829~$! ).
Thus we had to stop an initial systematic search for increasing values of $%
b,m,S_{0}$ : however eventually we tested the convergence for a huge number
of random chosen values of $b,m,S_{0}$.

Of course hope is hope, not reality: realistically one should expect to find
easily a large number of counter-examples, namely a large number of limit
cycles different from the Collatz type one given in (\ref{LC}) (divergences
also could possibly arise, but we never found one).

Indeed for $m=1$ it is easy to find counter-examples even for small values
of $b,S_{0}:$

\begin{itemize}
\item $b=3,S_{0}=5$ reaches the limit cycle%
\begin{equation}
7,30,10,42,14,57,19,78,26,105,35,141,47,189,63,21,7
\end{equation}

\item $b=4,S_{0}=11$ reaches the limit cycle%
\begin{eqnarray}
&&23,116,29,148,37,188,47,236,59,296,74,372,93,468,  \nonumber \\
&&117,588,147,736,184,45,232,58,292,73,368,92,23
\end{eqnarray}

\item $b=6,S_{0}=7$ reaches the limit cycle%
\begin{eqnarray}
&&23,162,27,192,32,228,38,270,45,318,53,372,  \nonumber \\
&&62,438,73,516,85,606,101,708,118,828,138,23
\end{eqnarray}

\item $b=9,S_{0}=31$ reaches the limit cycle%
\begin{eqnarray}
&&35,351,39,396,44,441,49,495,55,558,62,621,  \nonumber \\
&&69,693,77,774,86,864,96,963,107,1071,119,  \nonumber \\
&&1197,133,1332,148,1485,165,1656,184,1845,  \nonumber \\
&&205,2052,228,2286,254,2547,283,2835,315,35.
\end{eqnarray}
\end{itemize}

Many other counter-examples can be found for $m=1$; note however that we
found none not only for $b=2$ (this was expected since that is the original
Collatz map) but also for other small values of $b$ as $b=5,b=7,b=8$ (and
this was quite unexpected).

Let us go to $m=2.$ The first counter-example we found is the following one:

\begin{itemize}
\item $b=2,S_{0}=23$ reaches the limit cycle%
\begin{equation}
37,188,94,47,236,118,59,296,148,74,37
\end{equation}
\end{itemize}

Again it is easy to find other counter-examples for $b=2,m=2.\smallskip $

\textit{To our surprise, for }$m=2,b>2$\textit{\ and for }$m>2,b\geq 2$%
\textit{\ \ we have found no other counter-example}... \smallskip 

Since our computer investigation, even if massive, has not been systematic
(nor obviously exhaustive), we do not dare to propose a new general
conjecture of the Collatz type (however we dare to offer $\texteuro \ 1$ for
each of the first 100 counter-examples in the cases $%
m=1,b=5,7,8;m=2,b>2;m>2,b>=2$ ).

Nevertheless, motivated by some considerations about the behavior of some CA
(Cellular Automata) which mimic our generalized Collatz map (for CA related
to the original Collatz map see \cite{ca}), we dare to propose the following
conjecture (exactly analogous to the original Collatz conjecture):\smallskip

\begin{conjecture}
If $m=b-1$ the sequence (\ref{SeqB}) shall eventually reach the number $1$,
irrespective of which positive integer is chosen initially.\smallskip
\end{conjecture}

\textbf{Remark}\textit{\ \ \ Note that\bigskip\ for }$b=2$\textit{\ this is
just the original Collatz conjecture.}

We have tested this conjecture for a large number of values of $b$ and $%
S_{0} $. \ 

We offer $\texteuro 1$ for a counter-example, $\texteuro \ 100$ for a proof.

\bigskip

\textbf{Acknowledgements \ }\textit{\ It is a pleasure to thank F. Calogero
for useful discussions and M. Pistella for valuable help in the computer
investigation.\bigskip }

\end{document}